\documentclass[11pt,notitlepage]{article}
\usepackage{amsmath}
\usepackage{amssymb}
\usepackage{graphicx}
\usepackage{authblk}

\newtheorem{thm}{Theorem}

\newtheorem{prop}{Proposition}

\title{Stable Fixed Points of  Card Trick Functions}
\author[1]{Jyoti Champanerkar \thanks{\texttt{champanerkarj@wpunj.edu}}}
\author[2]{Mahendra Jani \thanks{\texttt{janim@wpunj.edu}}}
\affil[1,2]{Department of Mathematics, William Paterson University, Wayne, NJ 07470, USA.}
\date{}
\begin{document}
\maketitle
\abstract{
\noindent
The \textit{$21$-card trick} is a way of dealing cards in order to predict the card selected by a volunteer. We give a mathematical explanation of why the well-known 21-card trick works using a simple linear discrete function. The function has a stable fixed point which corresponds to the position where the selected card reaches at the end of the trick. We then generalize the $21(7 \times 3)$-card trick to a $p \times q$ - card trick where $p$ and $q$ are odd integers greater than or equal to three, determine the fixed point and prove that it is also stable.
}

\section{Introduction}
Mathematics of card tricks has been a fascinating and fun subject for a long time. 
 Ideas from number theory, group theory, dynamical systems and computer science have often been used to explain card tricks (\cite{Ledet}, \cite{MedMorr}, \cite{RamnathScully} and \cite{Scully}).  Conversely, playing cards have also been often used to illustrate mathematical concepts of probability distributions (\cite{Quinn} and \cite{VanHecke}) and group theory (\cite{Ensley} and \cite{Lynch}). The \textit{$21$-card trick} involves dealing the cards and picking them up repetitively till the selected card is located.\\ 
\ \\
The $21$-card trick has been studied by Amir-Moez(\cite{Amir}) and Bolker (\cite{Bolker}). While, Amir-Moez explains the trick as a limit of a function, there is no closed-form formula given for the function. 
 Bolker (\cite{Bolker}), explains the trick (referred to as Gergonne's trick in the paper) using ternary expansion and Radix sorting algorithm. \\
\ \\
In this paper, we give a mathematical explanation of why the $21$-card trick works using stability ideas from dynamical systems.  The trick is represented using a simple linear discrete function  which has a stable fixed point. It is shown that within three iterations any selected card reaches the position corresponding to the stable fixed point of the function. Further, two generalizations of the $21$-card trick are discussed, namely, the \textit{$p\times 3$-card trick} and the \textit{$p \times q$-card trick}. A functional representation is provided and existence of a stable fixed point is proved for each of the generalizations.\\
\ \\
\noindent
Since the distribution (laying down) and collection of cards are extremely important, we describe this process in section 2. In Section 3,  we define a linear discrete function for the 21-card trick, prove that it has a stable fixed point and that it is reached within three iterations irrespective of the initial position of the selected card.   In section 4, we generalize the $21$-card trick to $p \times 3$-card trick where $p$ is any odd integer greater than $1$. We determine the fixed point for the $p \times 3$ case and show that it can be reached within $a$ iterations (where $p=2a+1$), irrespective of the initial position of the selected card, using induction on $a$. In Section 5, we further generalize the $21$-card trick to $p \times q$-card trick where $p$ and $q$ are odd integers greater than $1$. 
The exact number of iterations required for the $p\times q$ case depends on how large $p$ is compared to $q$.\\
\ \\
\noindent
The proof in Section 3 is direct and easily constructible. This is generalized in Section 4, where the proof is more intricate and requires careful handling of inequalities. Finally, in Section 5 the method of induction is abstract and distinct from the one applied in Section 4. Further, this approach is pedagogically more appealing for undergraduates and shows them how mathematics is developed. Part of this work was presented at a Math Fair for high school students visiting a university campus, and was well-received.

\section{Method for distributing the cards}\label{sec:Method}
Development and properties of the card-trick functions described in this paper are based on how the cards are distributed each time. Although there is no sleight-of-hand involved, the order in which the cards are laid out and collected is important. The main idea is to place the column indicated by the volunteer, in the \textit{middle} of the stack before the next distribution. The method is described completely in this section and is referred to several times in the paper.\\
 \ \\
\noindent \textit{Hold $21$ cards face-up in one hand. Place the first $3$ cards face-up in a row from left to right. Then place the remaining cards in subsequent rows $2, 3, \cdots, 7$.
A volunteer is then asked to mentally select a card and  point out only the \textbf{column} in which the selected card lies. Suppose the volunteer points to column $3$ as the column  containing the selected card. Now collect all cards in column $1$ in a stack, with the bottom card as first and without disturbing the order. Then place cards from column $3$ on top of the stack, again picking the entire column at once, without disturbing the order. Now collect the cards in column $2$, and place them on top of the stack.
Redistribute as before, one card at a time, left to right, $3$ cards in each row. Ask the volunteer to point the column containing the previously marked card. Then stack as before with the column pointed out by the volunteer in the middle of the pile. }\\
\ \\
\noindent
The magic of the card trick is that within three iterations, the marked card will end up in the middle of the entire stack at position $11$. For a more detailed description and pictures about the method for distributing the cards, see \textit{Find the Card} (Magic Exhibit \cite{MagicExhibit}) or \textit{Do a $21$ card trick} (WikkiHow \cite{Wikki}).
\begin{figure}[!hbt]
\begin{center}

\begin{minipage}{0.2\textwidth}
\includegraphics[width=0.65\textwidth]{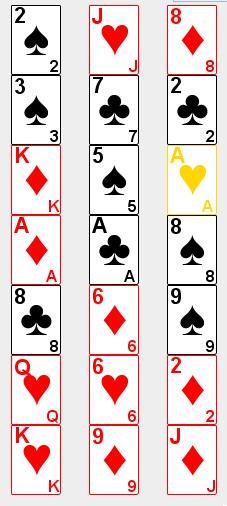}
\end{minipage}
\begin{minipage}{0.2\textwidth}
\includegraphics[width=0.65\textwidth]{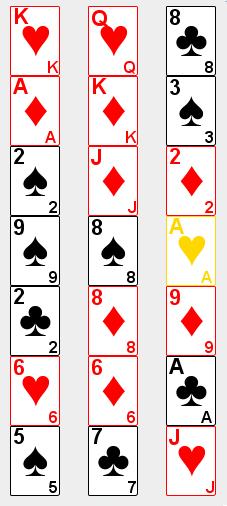}
\end{minipage}
\begin{minipage}{0.2\textwidth}
\includegraphics[width=0.65\textwidth]{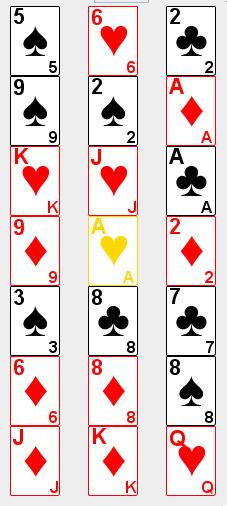}
\end{minipage}
\begin{minipage}{0.2\textwidth}
\includegraphics[width=0.5\textwidth]{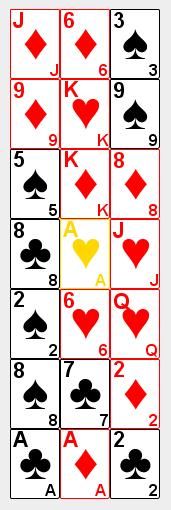}
\end{minipage}

\parbox{0.82\textwidth}{
\caption{\small{\textit{An Example: The selected card is \textit{Ace of Hearts} located in column 3 in the initial distribution. After the cards are collected with column 3 in the middle, and redistributed as described, \textit{Ace of hearts} is still in column 3. When the cards are collected with column 3 in the middle, and redistributed, \textit{Ace of Hearts} is in columns 2. Finally, the cards are collected with column 2 in the middle and redistributed, and \textit{Ace of Hearts} is the 11th card in the pile of 21 cards. In this case, the selected card reached the mid point in two iterations and stayed there.
\label{fig:fig01}}}
}
}
\end{center}
\end{figure}

\section{Function describing the $21$-card trick ($p=7$ and $q=3$)}
In this Section, we explain why the $21$-card trick works by using a mathematical function. For an introduction to discrete dynamical systems, fixed points, and stability properties please see Marotto (\cite{Marotto}).\\
\ \\
\noindent
Let $21$ cards be placed, face-up, in a $7$ rows $\times 3$ columns distribution. Let $n$ denote the position of the selected card. Each integer $n$, such that $1\leq n \leq 21$, can be expressed in the form $n=3k-l$ where $k$ and $l$  are unique integers with $1\leq k \leq 7$ and $0 \leq l \leq 2$. In the expression $3k-l$, integer $k$ indicates the row and integer $l$ is such that $3-l$ is the column number in which the card lies. That is, when $l=2$, the card is in the first column, when $l=1$ the card is in the second column and when $l=0$, the card is in the third column. Let $f(n)$ be the position of the selected card after the column containing the selected card is stacked in the middle of the pile, in the manner described in Section \ref{sec:Method}. Suppose that the selected card lies in row $k$ and in any column. When the cards are collected, the column containing the selected card becomes the middle column. The selected card is placed in the $k$th row of the middle column in the collected pile of cards. In the collected pile, there are $7$ cards from the first column, and $7-k$ cards from the middle column before the selected card, making its position $7+7-k+1$ or $15-k$. Hence, we define,
\begin{equation}
\boxed{f(n)= 15-k},
\label{eq:f_21a}\end{equation}
where $n=3k-l$, $1 \leq n \leq 21$ and $1\leq k \leq 7$.

\begin{thm}:
Let $n_0=3k_0 - l_0$ be the position of the selected card in the initial layout, where $1\leq n_0 \leq 21$,  $1 \leq k_0 \leq 7$ and $0 \leq l_0\leq 2$.  Then $f^3(n_0) = 11$, $\forall n_0$, where $f$ is as described in equation (\ref{eq:f_21a}). That is, $11$ is the stable fixed point of $f$. 
\end{thm}
\textbf{Proof: } We first show that $11$ is the fixed point of the function. Since $11=3(4)-1$, we have $k=4$. Substituting $n=11$ in equation (\ref{eq:f_21a}), we have $f(11)=15-4=11$. Thus, $11$ is a fixed point of the function $f$.\\
\ \\
\noindent We now prove the stability of the fixed point.
Let  $n_0=3k_0 - l_0$ be the position of the selected card in the initial layout. Then, $1 \leq n_0 \leq 21$, $1\leq k_0 \leq 7$  and $0 \leq l_0 \leq 2$.\\
\ \\
\noindent
Let $n_1 = f(n_0)= 15-k_0$. 
Since, $1\leq k_0 \leq 7$ we have, $-7 \leq -k_0 \leq -1$, hence, $15-7 \leq 15-k_0 \leq 15-1$, or  $8 \leq n_1 \leq 14$. 
Using $n_1 = 3 k_1 - l_1$ gives $3 \leq k_1 \leq 5$. \\
\ \\
\noindent
Let $n_2 = f^2(n_0) =  f(n_1)=15-k_1$. 
Since, $3 \leq k_1 \leq 5$, we have, $10 \leq n_2 \leq 12$. 
Now writing $n_2 = 3 k_2 - l_2$ gives $4 \leq k_2 \leq 4$ or $k_2=4$.\\
\ \\
\noindent
Finally,  let $n_3 = f^3(n_0) = f(n_2)= 15-k_2$. Since $k_2=4$, we have $n_3=f^3(n_0) =11$, as required.\\
\ \\
\noindent
Thus, $11$ is the stable fixed point of the $21$-card trick function $f$ in equation (\ref{eq:f_21a}) and the selected card reaches the $11$th position within three iterations irrespective of the initial position.$\square$

\begin{table}[!htb]
\begin{center}
\resizebox{0.8\textwidth}{!}
{\begin{tabular}{|r|c c c|c c c|c c c | c c c| c c c | c c c| c c c|}
\hline
 position of & & & & & & & & & & & & & & & & & & & & &\\
 selected card & 1 & 2 & 3 & 4 & 5 & 6 & 7 & 8 & 9 & 10 & 11 & 12 & 13 & 14 & 15 & 16 & 17 & 18 & 19 & 20 & 21 \\
 $n_0$ & & & & & & & & & & & & & & & & & & & & &\\
 \hline
 & & & & & & & & & & & & & & & & & & & & &\\
 $n_1=f(n_0)$ & 14 & 14 & 14 & 13 & 13 & 13 & 12 & 12 & 12 & 11 & 11 & 11 & 10 & 10 & 10 & 9 & 9 & 9 & 8 & 8 & 8\\
 & & & & & & & & & & & & & & & & & & & & &\\
\hline
& & & & & & & & & & & & & & & & & & & & &\\
$n_2=f^2(n_0)$ & 10 & 10 & 10 & 10 & 10 & 10 & 11 & 11 & 11 & 11 & 11 & 11 & 11 & 11 & 11 & 12 & 12 & 12 & 12 & 12 & 12\\
& & & & & & & & & & & & & & & & & & & & &\\
\hline
& & & & & & & & & & & & & & & & & & & & &\\
$n_3=f^3(n_0)$ & 11 & 11 & 11 & 11 & 11 & 11 & 11 & 11 & 11 & 11 & 11 & 11 & 11 & 11 & 11 & 11 & 11 & 11 & 11 & 11 & 11\\
& & & & & & & & & & & & & & & & & & & & &\\
\hline
\end{tabular}
}

\parbox{0.8\textwidth}
{
\caption{\small{\textit{This table enumerates where a selected card gets placed after each iteration, and how after three iterations the selected card reaches middle or the $11$th position for the $21$-card trick.}}
\label{table:selectedcards}}
}
\end{center}
\end{table}

\begin{figure}[!h]
\begin{center}
\includegraphics[width=0.6\textwidth]{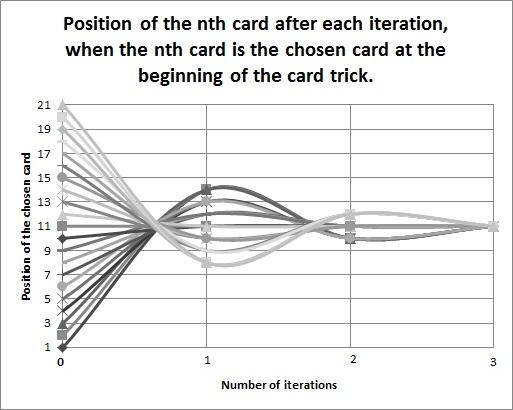}
\parbox{0.8\textwidth}
{\caption{\small{\textit{This figure is a graphical representation of Table \ref{table:selectedcards}}}.
}
}
\label{fig:stable_11}
\end{center}
\end{figure}

\newpage
\section{Generalization of $21$-card trick to $p \times 3$-card trick}
We extend the above result to $p \times 3$ cards where $p \geq 3$ is an odd integer. Note that, if $p=1$ then there is only one row of $3$ cards and the trick is redundant.

\begin{thm}: For an odd integer $p=2a+1$, with $a \geq 1$, suppose that $3p$ cards are distributed in a $p\times 3$ arrangement. Suppose that the selected card is at position $n$ and cards are collected and distributed as described in Section \ref{sec:Method}. Then within $a$ iterations, the selected card will be at position $\frac{3p+1}{2}$.
\end{thm}
\textbf{Proof: }
Let $p \times 3$ cards be placed, face-up, in $p$ rows and $3$ columns. Let $n$ denote the position of the selected card. Each integer $n$, such that $1\leq n \leq 3p$, can be expressed in the form
\begin{displaymath}n=3k-l\end{displaymath}
 where $1\leq k \leq p$ and $0 \leq l \leq 2$. In the expression $3k-l$, integer $k$ represents the row in which the card lies and integer $l$ is such that $3-l$ is the column number in which the card lies. That is, when $l=2$, the card is in the first column, when $l=1$ the card is in the second column and when $l=0$, the card is in the third column. Let $g(n)$ be the position of the selected card after the first iteration, where the cards are collected and redistributed as described in Section \ref{sec:Method}. Just as in the $21$-card trick, if the selected card is in row $k$, and in any column, when the cards are collected, it gets placed at row $k$ in the middle column in the stack. Thus, it is at position $p +p -k +1$ or $2p+1-k$ in the collected pile of cards. We define the function $g(n)$ as:
 \begin{equation}
\boxed{g(n)= 2p + 1 -k}.
\label{eq:g_3pa}\end{equation}

\noindent We first show that the midpoint $\frac{3p+1}{2}$ is a fixed point. Since
\begin{displaymath}
\frac{3p+1}{2} = \left(\frac{3p+3-3+1}{2}\right)=3\left(\frac{p+1}{2}\right) - 1,
\end{displaymath}
we have $k=\frac{p+1}{2}$ and $l=1$. Substituting $n=\frac{3p+1}{2}$ in equation (\ref{eq:g_3pa}), gives,
\begin{eqnarray*}
g\left(\frac{3p+1}{2}\right) &=& 2p + 1 -\frac{p+1}{2}\\
&=& \frac{4p+2-p-1}{2}\\
&=&\frac{3p+1}{2}
\end{eqnarray*}
Thus, the integer $\frac{3p+1}{2}$ is a \textit{fixed point} of the function $g$.\\
\ \\
\noindent
We now prove the stability of the fixed point. Let $n_0=3k_0-l_0$ be the initial position of the selected card. Then, $1 \leq n_0 \leq 3p$, $1 \leq k_0 \leq p$ and $0 \leq l_0 \leq 2$.\\
\ \\
\noindent
Let $n_1=g(n_0)=2p + 1 -k_0$. Since $1 \leq k_0 \leq p$, we have, $-p \leq -k_0 \leq -1$ and hence $2p + 1 - p \leq n_1 \leq 2p + 1 - 1$. Thus,
\begin{equation}\boxed{p+1 \leq n_1 \leq 2p}.\label{eq:3p_n1}\end{equation}
We find the bounds on $k_1$ by expressing $n_1=3k_1 -l_1$, in inequality (\ref{eq:3p_n1}),
\begin{eqnarray*}
p+1 &\leq 3k_1-l_1 \leq& 2p,\\
\text{thus, } p+1 +l_1 \leq 3k_1 &\text{ and }& 3k_1 \leq 2p+l_1.
\end{eqnarray*}
Since $0 \leq l_1 \leq 2$, we have $p+1 \leq 3k_1$ and $3k_1 \leq 2p+2$. Consider the first expression $p+1 \leq 3k_1$. This implies that $\frac{p+1}{3} \leq k_1$. Since, $p \geq 3$, $\frac{p+1}{3} \geq 1.33$. But $k_1$ is an integer, hence we get the lower bound $2\leq k_1$. Now consider the expression, $3k_1 \leq 2p+2$, which implies $k_1 \leq \frac{2}{3}(p+1)$. Given that $p\geq 3$, we consider the following two cases. (i) When $p=3$, $\frac{2}{3}(p+1)=8/3$ and since $k_1$ is an integer less than $\frac{2}{3}(p+1)$, it must be less than or equal to $2$, and $2=p-1$. (ii) Otherwise, $\frac{2}{3}(p+1)\leq p-1$ if and only if $5\leq p$. Thus we have,
\begin{equation} \boxed{2 \leq k_1 \leq p-1}.\label{eq:3p_k1}\end{equation}
Let $n_2=g(n_1)=2p + 1 -k_1$. If $p=3$, then $2\leq k_1 \leq p-1$ implies that $k_1=2$ and hence $n_2=5$. Thus proving that $\frac{3p+1}{2}$ is a stable fixed point of $g$ and that for all integers $n_0$, such that $1\leq n_0 \leq 9$, the fixed point is reached in two iterations.\\
\ \\
For the rest of the proof we assume $p\geq 5$ is an odd integer. As before using the bounds on $k_1$, we obtain the estimate \boxed{p+2\leq n_2 \leq 2p-1} and express $n_2$ as $n_2=3k_2 - l_2$. This implies that,
\begin{eqnarray*}
p+2 &\leq 3k_2-l_2 \leq& 2p-1\\
p+2 \leq 3k_2-l_2 \leq 3k_2 && 3k_2 - 2 \leq 3k_2 - l_2 \leq 2p-1\\
\text{or } \frac{p+2}{3}  \leq k_2 & \text{ and }& k_2 \leq \frac{2p+1}{3}
\end{eqnarray*}
Since, $p\geq 5$, we have $k_2 \geq 7/3$. But $k_2$ is an integer, hence, $k_2 \geq 3$. Also $\frac{2p+1}{3} \leq p-2$ if and only if $7 \leq p$. When $p=5$, we have $\frac{2p+1}{3}=11/3$ and since $k_2$ is an integer, $k_2 \leq 3 = p-2$. Hence we have the estimate
\begin{equation}\boxed{3 \leq k_2 \leq p-2}.\label{eq:3p_k2}\end{equation}
\noindent
Let $n_3=g(n_2)=2p + 1 -k_2$. If $p=5$, then $3\leq k_2 \leq p-2$ implies that $k_2=3$ and hence $n_3=8$. Thus proving that $\frac{3p+1}{2}$ is a stable fixed point of $g$ and that for all integers $n_0$, such that $1\leq n_0 \leq 15$, the fixed point is reached in three iterations.\\
\ \\
Using the estimate for $k_2$, we obtain \boxed{p+3 \leq n_3 \leq 2p-2}. Using $n_3=3k_3 - l_3$ and reasoning as before one can observe that \boxed{4 \leq k_3 \leq p-3}. If $p=7$, we obtain $k_3 =4$ and $n_4 =11$ and hence the claim holds for $p=7$. Results obtained for odd integers $3, 5, 7$ are summarized below in Table \ref{tab:computations_3p}.
\begin{table}[!hbt]
\begin{center}
\begin{tabular}{|c|c|c|c|c|}
\hline
 & Odd integer & fixed point & penultimate row &  last iteration \\
$a$ & $p=2a-1$ & $\frac{3p+1}{2}$ & $k_{a-1}$ & $n_a=2p+1-k_{a-1}$\\
\hline
$2$ &  $3$ & $5$ & $2$ & $5$\\
\hline
$3$ & $5$ & $8$ & $3$ & $8$\\
\hline
$4$ & $7$ & $11$ & $4$ & $11$\\
\hline
\end{tabular}
\caption{\small{\textit{Last iteration and fixed point computations for first few odd numbers $p$.}}\label{tab:computations_3p}}
\end{center}
\end{table}

\noindent
We will prove the theorem by induction. We have shown that for $p=3$, $2 \leq k_1 \leq p-1$. Assume that for any odd integer  $p \leq 2a-1$, with $a \geq 2$, the following inequality holds: \boxed{a \leq k_{a-1} \leq p-a+1} and we prove that a corresponding inequality holds for $p=2a+1$.\\
\ \\
\noindent
By definition, $n_a=2p+1 - k_{a-1}$. By induction hypothesis, if $p=2a-1$, then $k_{a-1}=a$ and $n_a=3a-1 = \frac{3p+1}{2}$. \\ 
\ \\
\noindent
Using the estimates for $k_{a-1}$, we have $p+a \leq n_a \leq 2p+1 -a$. Expressing $n_a = 3k_a - l_a$, yields $p+a \leq 3k_a$ and $3k_a \leq 2p+3 -a$. The inequality $a + 1/3 \leq k_a$ together with the fact that $k_a$ is an integer, implies that $a+1\leq k_a$. Now consider $ k_a \leq \frac{2p+3 - a}{3}$. Observe that $\frac{2p+3 - a}{3} \leq p-a$ whenever $2a+3 \leq p$. And if $p=2a+1$, then $\frac{2p+3 - a}{3}=a + 5/3$. Again using that $k_a$ is an integer, we have
$k_a \leq a+1 = p-a$. Thus, for all odd integers $p$ such that $p \geq 2a + 1$, we have the bounds
\begin{equation}\boxed{a +1 \leq k_a \leq p-a}.\end{equation}
In particular, when $p=2a+1$, $k_a = a+1$ and $n_a = 3a+2 = \frac{3p+1}{2}$. Thus for any $n_0$, $1 \leq n_0 \leq 3p$, $g^a(n_0)=\frac{3p+1}{2}$. Hence the fixed point $\frac{3p+1}{2}$ is stable.$\square$\\
\ \\
\noindent
In the $p \times 3$ case, it is observed that when $p=2a+1$, at most $a$ iterations are needed. In particular, when $p=7$ as in the $21$-card trick, $a=3$, and at most $3$ iterations are needed. See Table \ref{table:selectedcards} and Figure \ref{fig:stable_11}.

\section{Generalization of $21$-card trick to $p \times q$-card trick}
We now extend the above results to $p \times q$ cards where $p$ and $q$ are odd integers, greater than or equal to $3$.\\
\ \\
\noindent
Let $p \times q$ cards be be placed, face up as described in the Section \ref{sec:Method}. Let $n$ denote the position of the selected card. Each integer $n$, such that $1 \leq n \leq pq$, can be expressed in the form
\begin{displaymath}
n=qk - l
\end{displaymath}
where $1 \leq k \leq p$ and $0 \leq \l \leq q-1$. In the expression $qk-l$, integer $k$ represents the row in which the card lies. The integer $l$ is such that $q-l$ is the column number in which the card lies. That is, when $l=q-1$, the card is in the first column, when $l=q-2$ the card  in the second column etc., and when $l=0$, the card is in the $q$th column.\\
\ \\
\noindent Let $h(n)$ be the position of the selected card after the first iteration. Just as in the $21$-card trick, when the indicated column is placed in the middle of the collected stack of cards, there are $(q-1)/2$ columns of $p$ cards each, on top, then the indicated column, and then $(q-1)/2$ columns of $p$ cards each, below. If the selected card was in row $k$, in any column, it gets placed at position $p(q-1)/2+ p-k +1$ in the collected stack of cards. Hence, we define $h(n)$ as:
\begin{equation}
\boxed{h(n)=\left(\frac{q+1}{2}\right)p + 1 - k}.
\label{eq:h_pqa}
\end{equation}
\begin{thm}:
For $n=qk-l$  and $h(n)$ as defined in equation (\ref{eq:h_pqa}), the midpoint, $\frac{pq+1}{2}$ is a fixed point of $h$.
\end{thm}
\noindent \textbf{Proof:} We first observe that
\begin{displaymath} \frac{pq+1}{2}= q \left(\frac{p+1}{2}\right) - \left(\frac{q-1}{2}\right),
\end{displaymath}
giving $k=\frac{p+1}{2}$. Substituting $n=\frac{pq+1}{2}$  and $k=\frac{p+1}{2}$ in equation (\ref{eq:h_pqa}), gives,
\begin{eqnarray*}
h\left(\frac{pq+1}{2}\right) &=& \left(\frac{q+1}{2}\right)p + 1 - \frac{p+1}{2}\\
&=&\frac{pq+p+2-p-1}{2}\\
&=& \frac{pq+1}{2}.
\end{eqnarray*}
Thus, the integer $\frac{pq+1}{2}$ is a fixed point of $h$.$\square$

\subsection{Stability of the fixed point}
 We prove that the fixed point $\frac{pq+1}{2}$ of $h$ is in fact stable using several propositions proved in this subsection.

\begin{prop}:
Let $h(n)$ be the function as described in equation (\ref{eq:h_pqa}). If
\begin{displaymath}
q\left(\frac{p+1}{2}\right) - (q-1) \leq n \leq q\left(\frac{p+1}{2}\right),
\end{displaymath}
then
\begin{displaymath}
h(n) = \frac{pq+1}{2}.
\end{displaymath}\label{prop:pqmiddlerow}
\end{prop}
\noindent \textbf{Proof:} For $n$ as given, $k=\frac{p+1}{2}$, substituting in $h(n)$ in equation (\ref{eq:h_pqa}) gives,
\begin{eqnarray*}
h(n)&=& \left(\frac{q+1}{2}\right)p + 1 - \frac{p+1}{2}\\
&=& \frac{pq+1}{2}
\end{eqnarray*}
as claimed. $\square$\\
\ \\
\noindent
Proposition \ref{prop:pqmiddlerow}, implies that, if at any time the position of the card is in the middle row, that is the $\frac{p+1}{2}$th row, then in the next iteration, the card reaches the fixed point irrespective of the column it was in, in the middle row.

\begin{prop}:
Let $h(n)$ be the function as described in equation (\ref{eq:h_pqa}). If $n=qk-l$ is such that,
\begin{displaymath}
\frac{p+1}{2} - \frac{q-1}{2} \leq k \leq \frac{p+1}{2} + \frac{q-1}{2},
\end{displaymath}
then
\begin{displaymath}
h^2(n)=\frac{pq+1}{2}.
\end{displaymath}\label{prop:pqmiddleblock}
\end{prop}
\noindent \textbf{Proof:} For $n=qk -l$ in rows $\frac{p+1}{2} - \frac{q-1}{2}$ through $\frac{p+1}{2} + \frac{q-1}{2}$, we have,
\begin{eqnarray*}
\frac{p+1}{2} - \frac{q-1}{2} \leq &k& \leq \frac{p+1}{2} + \frac{q-1}{2}\\
\text{ hence, } \left(\frac{q+1}{2}\right)p + 1 - \left(\frac{p+1}{2} + \frac{q-1}{2}\right) \leq &h(n)& \leq \left(\frac{q+1}{2}\right)p + 1 - \left(\frac{p+1}{2} - \frac{q-1}{2}\right)\\
\text{ which simplifies to, } \frac{pq+p+2-p-1-q+1}{2} \leq &h(n)& \leq \frac{pq+p+2-p-1+q-1}{2}\\
\frac{pq+2-q}{2} \leq &h(n)& \leq \frac{pq+q}{2}\\
q\left(\frac{p+1}{2}\right)-(q-1) \leq &h(n)& \leq q\left(\frac{p+1}{2}\right).
\end{eqnarray*}
That is $h(n)$ lies in the middle row or the $\frac{p+1}{2}$th row. By Proposition \ref{prop:pqmiddlerow}, $h^2(n)=h(h(n))=\frac{pq+1}{2}$.$\square$

\subsubsection{Remarks}
\begin{enumerate}
\item Proposition \ref{prop:pqmiddleblock} implies that, if after any iteration, the position of the selected card is within the block of $q$ rows containing the middle row, that is within the $\frac{q-1}{2}$ rows above or below the middle row and the middle row itself, then the card reaches the fixed point within two more iterations, irrespective of the column it was in.

\item If $p\leq q$, then all $p$ rows are within the middle block of $q$ rows, and within two iterations, the selected card will be at the position denoted by the fixed point $\frac{pq+1}{2}$. Hence, $\frac{pq+1}{2}$ is a \textit{stable} fixed point for the function $h$ for the case when $p \leq q$.
\end{enumerate}

\begin{figure}[!hbt]
 \begin{center}
 \includegraphics[width=0.7\textwidth]{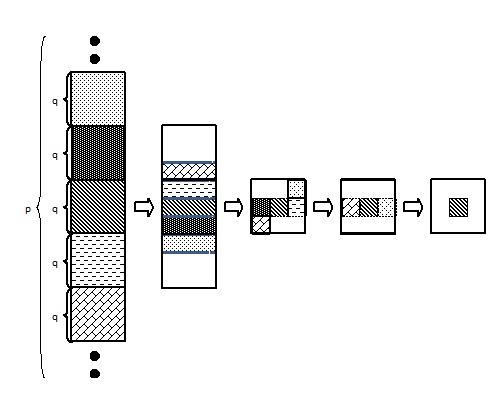}
\parbox{0.8\textwidth}{
 \caption{\small{\textit{Pictorial representation of how blocks of $q$ rows  starting in the middle, get mapped to a single row. The middle block of $q$ rows gets mapped to the middle row via function $h$ and the middle row to the middle position.}}}
 \label{fig:qblockdiagram}
 }
 \end{center}
 \end{figure}

 \noindent In general, for $p>q$, we think of $p$ rows to be divided into blocks of $q$ rows, starting at the middle row and extending to the first and last rows. Thus, the middle block consists of the middle row ($\frac{p+1}{2}th$ row), and $\frac{q-1}{2}$ rows above and $\frac{q-1}{2}$ rows below the middle row. Since $p-q$ is an even number, there are an equal number of rows above and below the middle block of $q$ rows. Hence there are an equal number of blocks above and below the middle block with the possibility that the blocks containing the first and last row may have less than or equal to $q$ rows. We have shown in the previous proposition that if the selected card lies in the middle block of $q$ rows, then within two iterations, it would have the reached the fixed point position $\frac{pq+1}{2}$. We will now show by induction that if the selected card lies outside the middle block then it reaches the fixed point position $\frac{pq+1}{2}$.

 \begin{thm}: Let $p, q \geq 3$ be odd integers. For $pq$ cards distributed as described in Section \ref{sec:Method}, let $n=qk-l$ denote the position of the selected card where as before, $1 \leq n \leq pq$, $1\leq k \leq p$ and $0\leq l \leq q-1$. Let $h(n)$ be as described in equation (\ref{eq:h_pqa}). For $k$ such that,
 \begin{eqnarray*}
 \frac{p+1}{2} + \frac{q-1}{2} + jq + 1 \leq &k& \leq \frac{p+1}{2} + \frac{q-1}{2} + jq + q,
 \end{eqnarray*}
 where $j=0,1,2 \cdots$ then $h\circ h \cdots h(n)=\frac{pq+1}{2}$ in a finite number of compositions. That is, $\frac{pq+1}{2}$ is a stable fixed point of $h$.
 \label{thm:stable_pq}
 \end{thm}
The bounds for $k$ mean that the selected card initially lies within the block consisting of $jq +1$ and $jq+q$ rows below the middle block.\\
\ \\
\noindent \textbf{Proof: } Let $\mathcal{P}(m)$ be the proposition that:
 If \begin{eqnarray*}
 \frac{p+1}{2} + \frac{q-1}{2} + mq + 1 \leq &k& \leq \frac{p+1}{2} + \frac{q-1}{2} + mq + q,
 \end{eqnarray*}
 then $h\circ h \cdots h(n)=\frac{pq+1}{2}$ in a finite number of compositions.\\
 The inductive hypotheses is that: $\mathcal{P}(m)$ is true for all $m$ such that $0 \leq m < i$.\\ We will show that $\mathcal{P}(i)$ is true.\\
 Let $n=qk-l$. Let $k$ be such that,
 \begin{eqnarray*}
 \frac{p+1}{2} + \frac{q-1}{2} + iq + 1 \leq &k& \leq \frac{p+1}{2} + \frac{q-1}{2} + iq + q.
 \end{eqnarray*}
 Then,
 \begin{eqnarray*}
 \frac{p+q}{2} +  iq + 1 \leq &k& \leq \frac{p+q}{2} + iq + q\\
 -\left(\frac{p+q}{2} + iq + q \right) \leq & -k & \leq -\left(\frac{p+q}{2} + iq + 1\right)\\
 \left(\frac{q+1}{2}\right)p + 1 -\left(\frac{p+q}{2} + iq + q \right) \leq & h(n) & \leq \left(\frac{q+1}{2}\right)p + 1 -\left(\frac{p+q}{2} + iq + 1\right)\\
 \frac{pq+p+2}{2} - \frac{p+q+2iq + 2q}{2} \leq & h(n)& \leq \frac{pq+p+2}{2} - \frac{p+q+2iq + 2}{2}\\
 \frac{pq-q-2iq}{2} + 1 -q \leq & h(n) & \leq  \frac{pq-q-2iq}{2} \\
 \left( \frac{p+1}{2} - (i+1)\right)q - (q-1) \leq & h(n) & \leq \left( \frac{p+1}{2} - (i+1)\right)q.
 \end{eqnarray*}
Thus $h(n)$ lies in the row $\frac{p+1}{2} - (i+1)$. That is, if the selected card was anywhere in the entire block of $q$ rows, it will be placed at row $\frac{p+1}{2} - (i+1)$ in the next iteration. Since, $\frac{p+1}{2} - (i+1) < \frac{p+1}{2} + \frac{q-1}{2} + iq + 1$, by induction hypothesis, the card will reach the position denoted by the fixed point in a finite number of iterations. By the principle of induction, $\mathcal{P}(j)$ is true for all $j=0,1,2 \cdots$. Hence, $\frac{pq+1}{2}$ is a stable fixed point of the function $h$, for all odd $p,q \geq 3$.$\square$\\
\ \\
\noindent
We proved in Proposition \ref{prop:pqmiddleblock}, that for $p \leq q$, at most two iterations are needed for any selected card to reach the middle of the stack. For $p > q$, we were able to prove in Section 5, that a finite number of iterations are required, however, the exact number of iterations depends on how large $p$ is compared to $q$. We expect that, for $q+1 \leq p \leq 3q$ three iterations are required, and for $3q+1 \leq p \leq 5q$ four iterations are required, for any selected card to reach the middle of the stack.

\section*{Acknowledgements}
 Thanks to Dr. M. Zeleke for suggestions and modifications; Work on this project was partially supported by Assigned Release Time (ART), William Paterson University awarded to both authors.

\end{document}